\label{key}
\documentclass[12pt, reqno]{amsart}
\usepackage{amssymb, amsthm, amsmath, amsfonts}
\usepackage{array, epsfig}
\usepackage{bbm}
\usepackage{hyperref}
\usepackage[numbers,square]{natbib}
\usepackage{color}

  \evensidemargin
\oddsidemargin
\newtheorem{theorem}{Theorem}[section]

\newtheorem{corollary}[theorem]{Corollary}
\newtheorem{proposition}[theorem]{Proposition}

\newcommand{\E}{\mathbb{E}\,}

\newcommand{\R}{\mathbb{R}}
\newcommand{\Z}{\mathbb{Z}}

\newcommand{\bx}{\mathbf{x}}

\newcommand{\bp}{\mathbf{e}}

\newcommand{\by}{\mathbf{y}}

\newcommand{\dd}{{\rm d}}

\newcommand{\ind}{\mathbbm{1}}

\newcommand{\El}{\mathcal{E}}
\newcommand{\B}{\mathbb{B}}

\newcommand{\eqdistr}{\stackrel{d}{=}}
\newcommand{\conv}{\mathop{\mathrm{conv}}\nolimits}
\newcommand{\SPAN}{\mathop{\mathrm{span}}\nolimits}

\newcommand{\tc}{\textcolor{red}}

\title{Random affine simplexes}

\thanks{The work was done with the financial support of the Bielefeld University (Germany). The work of F.G. and D.Z. is supported by the grant SFB 1283. The work of A.G. is supported by the grant IRTG 2235. The work of D.Z. is supported by the grant RFBR 16-01-00367 and by the Program of Fundamental Researches of Russian Academy of Sciences ``Modern Problems of Fundamental Mathematics''.}

\author{Friedrich G\"otze}
 \address{Friedrich G\"otze: Faculty of Mathematics,
 Bielefeld University,
 P. O. Box 10 01 31,
 33501 Bielefeld, Germany}
 \email{goetze@math.uni-bielefeld.de}

\author{Anna Gusakova}
\address{Anna Gusakova: Faculty of Mathematics,
 Bielefeld University,
 P. O. Box 10 01 31,
 33501 Bielefeld, Germany}
\email{agusakov@math.uni-bielefeld.de}

\author{Dmitry Zaporozhets}
\address{Dmitry Zaporozhets: St.\ Petersburg Department of Steklov Mathematical Institute,
Fontanka~27,
191023 St.\ Petersburg, Russia}
\email{zap1979@gmail.com}

\begin{document}

\begin{abstract}
For a fixed $k\in\{1,\dots,d\}$ consider arbitrary random vectors $X_0,\dots, X_{k}\in\mathbb R^d$ such that
the $(k+1)$-tuples $(UX_0,\dots, UX_{k})$  have the same distribution for any rotation $U$. Let $A$ be any non-singular $d\times d$ matrix. We show that the $k$-dimensional volume of the convex hull of affinely transformed $X_{i}$'s satisfies 
\[
|\conv(AX_0,\dots,AX_{k})|\eqdistr\frac{|P_\xi\El|}{\kappa_k}\cdot|\conv(X_0,\dots,X_{k})|,
\]
where 
$
\El:=\{\bx\in\R^d:{\bx^\top (A^\top A)^{-1}\bx}\leq 1\}
$
is an ellipsoid, $P_\xi$ denotes the orthogonal projection to a uniformly chosen random $k$-dimensional linear subspace $\xi$ independent of $X_0,\dots, X_{k}$, and $\kappa_k$ is the volume of the unit $k$-dimensional ball.

As an application, we derive the following integral geometry formula for ellipsoids:
\[
\frac{\kappa_{d}^{k+1}}{\kappa_k^{d+1}}\,\frac{\kappa_{k(d+p)+k}}{\kappa_{k(d+p)+d}}\,\int\limits_{A_{d,k}}|\El\cap E|^{p+d+1}\,\mu_{d,k}(dE)=|\El|^{k+1}\,\int\limits_{G_{d,k}}|P_L\El|^p\,\nu_{d,k}(dL),
\]
where $p>-1$ and $A_{d,k}$ and $G_{d,k}$ are the affine and the linear Grassmannians equipped with their respective  Haar measures. The case $p=0$ reduces to an affine version of the integral formula of Furstenberg and Tzkoni~\cite{FT71}.
\end{abstract}

\maketitle

\section{Main results}
\subsection{Basic notation}\label{1733}

First we introduce some basic notion of integral geometry following~\cite{SW08}. The Euclidean space $\R^d$ is equipped with the Euclidean scalar product $\langle \cdot,\cdot\rangle$. The volume is denoted by $|\cdot|$. Some of the sets we consider have dimension  less than $d$. In fact, we consider 3 classes: the convex hulls of $k+1$ points, orthogonal projections to $k$-dimensional linear subspaces, and intersections with $k$-dimensional affine subspaces, where $k\in\{0,\dots,d\}$. In this case,  $|\cdot|$ stands for the $k$-dimensional volume.

The unit ball  in $\R^k$ is denoted by $\B^k$. For $p>0$ we write 
\begin{equation}\label{1409}
\kappa_p:=\frac{\pi^{p/2}}{\Gamma\left(\frac p2+1\right)},
\end{equation}
where for an  integer $k$ we have $\kappa_k=|\B^k|$.

For $k\in\{0,\dots,d\}$, the linear (resp., affine) Grassmannian of $k$-dimensional linear (resp., affine) subspaces of $\R^d$ is denoted by $G_{d,k}$ (resp., $A_{d,k}$) and is equipped with a unique rotation invariant (resp., rigid motion invariant) Haar measure $\nu_{d,k}$ (resp., $\mu_{d,k}$), normalized by
\[
\nu_{d,k}(G_{d,k})=1,
\]
and
\[
\mu_{d,k}\left(\left\{E\in A_{d,k}:\,E\cap \B^d\ne\emptyset\right\}\right)=\kappa_{d-k},
\]
respectively.

A compact convex subset $K$ of $\R^d$ with non-empty interior is called a convex body. We define the intrinsic volumes of $K$  by Kubota's formula:
\begin{equation}\label{1309}
V_k(K)={d \choose k}\,\frac{\kappa_d}{\kappa_k\kappa_{d-k}}\,\int\limits_{G_{d,k}}|P_LK|\,\nu_{d,k}(dL),
\end{equation}
where $P_LK$ denotes the image of $K$ under the orthogonal projection to $L$.

For  $L\in G_{d,k}$ (resp., $E\in A_{d,k}$) we denote by $\lambda_L$ (resp., $\lambda_E$) the $k$-dimensional Lebesgue measures on $L$ (resp., $E$).

\subsection{Affine transformation of spherically symmetric distribution}\label{013}

For a fixed $k\in\{1,\dots,d\}$ consider random vectors $X_0,\dots, X_{k}\in\R^d$ (not necessarily i.i.d.) with an arbitrary \emph{spherically symmetric} joint distribution. By this we mean that the $(k+1)$-tuple $(UX_0,\dots, UX_{k})$  has the same distribution for any orthogonal $d\times d$ matrix $U$.
The  convex hull
\[
\conv(X_0,\dots,X_{k})
\]
is a $k$-dimensional simplex (maybe degenerate) with well-defined $k$-dimensional volume
\begin{equation}\label{1428}
|\conv(X_0,\dots,X_{k})|.
\end{equation}

How does the volume in~\eqref{1428} change under affine transformations? For $k=d$ the answer is obvious: it is multiplied by the  determinant of the transformation. The case $k<d$ presents a more delicate problem.
\begin{theorem}\label{1305}
Let $A$ be any non-singular $d\times d$ matrix and $\El$ be the ellipsoid defined by
\begin{equation}\label{2259}
\El:=\left\{\bx\in\R^d:{\bx^\top (A^\top A)^{-1}\bx}\leq 1\right\}.
\end{equation}
Then we have
\begin{equation}\label{1514}
|\conv(AX_0,\dots,AX_{k})|\eqdistr\frac{|P_\xi\El|}{\kappa_k}\cdot|\conv(X_0,\dots,X_{k})|,
\end{equation}
where $P_\xi$ denotes the orthogonal projection to a uniformly chosen random $k$-dimensional linear subspace $\xi$ independent of $X_0,\dots, X_{k}$.
\end{theorem}

Due to Kubota's formula (see~\eqref{1309}) $\E|P_\xi\El|$ is proportional to  $V_k(\El)$. Thus, taking ex\-pectation in~\eqref{1514} and using the formula
\[
V_k(\B^d)={d \choose k}\frac{\kappa_d}{\kappa_{d-k}}
\]
readily implies the following corollary.
\begin{corollary}\label{1721}
Under the assumptions of Theorem \ref{1305} we have
\begin{equation}\label{2314}
\E|\conv(AX_0,\dots,AX_{k})|=\frac{V_k(\El)}{V_k(\B^d)}\,\E|\conv(X_0,\dots,X_{k})|.
\end{equation}
\end{corollary}

For a formula of $V_k(\El)$ see~\cite{KZ14}. Relation~\eqref{2314} can be generalized to higher moments using the notion of \emph{generalized} intrinsic volumes introduced in~\cite{DP12}, but we shall skip to describe details here.

\bigskip

The main ingredient of the proof of Theorem~\ref{1305} is the following deterministic version of \eqref{1514}.

\begin{proposition}\label{1843}
Let $A$ and $\El$ be as in Theorem~\ref{1305}.
Consider $\bx_1,\dots,\bx_k\in\R^d$ and denote by $L$ their span (linear hull).  Then
\begin{equation}\label{2301}
	|\conv(0,A\bx_1,\ldots,A\bx_k)|=\frac{|P_L\El|}{\kappa_k}\cdot|\conv(0,\bx_1,\ldots,\bx_k)|.
\end{equation}	
\end{proposition}

Let us stress that here, the origin is added to the convex hull.

Applying~\eqref{2301} to standard Gaussian vectors (details are in Subsection~\ref{016}) leads to the following representation.
\begin{corollary}\label{009}
Under the assumptions of Theorem \ref{1305} we have
\begin{equation}\label{2306}
\frac{|P_\xi\El|}{\kappa_k}\eqdistr \left({\frac{\det\big(G^\top A^\top A G\big)}{\det\big(G^\top G\big)}}\right)^{1/2}\eqdistr \left({\frac{\det\big(G_\lambda^\top G_\lambda\big)}{\det\big(G^\top G\big)}}\right)^{1/2},
\end{equation}
where $G$ is a random $d\times k$ matrix with i.i.d. standard  Gaussian entries $N_{ij}$ and $G_\lambda$ is a random $d\times k$ matrix with the entries   $\lambda_i N_{ij}$, where  $\lambda_1,\dots,\lambda_d$ denote the singular values of $A$.
\end{corollary}
 Thus, we obtain the following version of~\eqref{1514}.
 \begin{corollary}
Under the assumptions of Theorem \ref{1305} and Corollary \ref{009}  we have
\begin{align*}
|\conv(AX_0,\dots,AX_{k})|&\eqdistr\left({\frac{\det\big(G^\top A^\top A G\big)}{\det\big(G^\top G\big)}}\right)^{1/2}|\conv(X_0,\dots,X_{k})|\\
&\eqdistr\left({\frac{\det\big(G_\lambda^\top G_\lambda\big)}{\det\big(G^\top G\big)}}\right)^{1/2}|\conv(X_0,\dots,X_{k})|.
\end{align*}
\end{corollary}
The important special case $k=1$ corresponds to the distance between two random points.
\begin{corollary}
Under the assumptions of Theorem \ref{1305} we have
\[
|AX_0-AX_1|\eqdistr\sqrt{\frac{\lambda_1^2N_1^2+\dots+\lambda_d^2N_d^2}{N_1^2+\dots+N_d^2}}\cdot|X_0-X_1|,
\]
where  $N_1,\dots,N_d$ are i.i.d. standard Gaussian variables and $\lambda_1,\dots,\lambda_d$ denote the singular values of $A$.
\end{corollary}

\subsection{Random points in ellipsoids}

Now suppose that  $X_0,\dots, X_{k}$ are independent and uniformly distributed in some convex  body  $K\subset\R^d$. A classical problem of stochastic geometry is to find the distribution of~\eqref{1428} starting with  its moments
\begin{equation}\label{1731}
\E|\conv(X_0,\dots,X_{k})|^p=\frac{1}{|K|^{k+1}}\,\int\limits_{K^{k+1}}{|\conv(\bx_0,\ldots,\bx_k)|^p\,\dd\bx_0\ldots \dd\bx_k}.
\end{equation}

The most studied case is $d=2, k=p=1$, when the problem reduces to the calculating the mean distance between two  uniformly chosen random points in a planar convex set  (see~\cite{eB25}, \cite{bG51}, \cite[Chapter~4]{lS76}, \cite[Chapter~2]{aM99}, \cite{uB14}).

For an arbitrary $d$ and $k=1$, there is an electromagnetic interpretation of~\eqref{1731} (see~\cite{HR04}): a transmitter $X_0$ and a receiver $X_1$ are placed uniformly at random in $K$. It is  empirically known that the  power received decreases with an inverse distance law of the form $1/|X_0-X_1|^\alpha$, where $\alpha$ is the so-called path-loss exponent, which depends on the environment in which both are located (see~\cite{RABT02}). Thus\tc{,} with $k=1$ and $p=-n\alpha$ Equation~~\eqref{1731} expresses the $n$-th moment of the  power received ($n<d/\alpha$).

The case of arbitrary $k$ and $d$ was studied only for $K$ being a ball. In~\cite{rM71} it was shown (see also~\cite[Theorem~8.2.3]{SW08}) that for  $X_0,\dots, X_{k}$ uniformly distributed in the unit ball $\B^d\subset\R^d$ and for an integer $p\ge 0$
\begin{equation}\label{018}
\E|\conv(X_0,\dots,X_{k})|^p=\frac{1}{(k!)^p}\,\frac{\kappa_{d+p}^{k+1}}{\kappa_{d}^{k+1}}\,\frac{\kappa_{k(d+p)+d}}{\kappa_{(k+1)(d+p)}}\,\frac{b_{d,k}}{b_{d+p,k}},
\end{equation}
where   $\kappa_k$ is defined in~\eqref{1409}  and for any real number $q>k-1$ we write (see~\cite[Eq.~(7.8)]{SW08})
\begin{equation}\label{1846}
b_{q,k}:=\frac{\omega_{q-k+1}\cdots\omega_{q}}{\omega_1\cdots\omega_k}
\end{equation}
with $\omega_p:=p\kappa_{p}$ being equal to the area of the unit $(p-1)$-dimensional sphere when $p$ is integer.

In~\cite[Proposition~2.8]{KTT17} this relation was extended to all real $p>-1$. It should be noted that Proposition~2.8 in \cite{KTT17} is formulated for real $p\ge 0$ only, but in the proof (see~p.~23) it is argued that by analytic continuation, the formula holds for all real $p>-1$ as well. Theorem~\ref{1305} implies (for details see Subsection~\ref{020}) the following generalization of~\eqref{018} for ellipsoids. Recall that $P_\xi$ denotes the orthogonal projection to a uniformly chosen random $k$-dimensional linear subspace $\xi$ independent of $X_0,\dots, X_{k}$.

\begin{theorem}\label{2227}
For  $X_0,\dots, X_{k}$ uniformly distributed in some non-degenerate ellipsoid $\El\subset\mathbb{R}^d$ and any real number $p>-1$ we have
\begin{equation}\label{1124}
\E|\conv(X_0,\dots,X_{k})|^p=\frac{1}{(k!)^p}\,\frac{\kappa_{d+p}^{k+1}}{\kappa_{d}^{k+1}}\,\frac{\kappa_{k(d+p)+d}}{\kappa_{(k+1)(d+p)}}\,\frac{b_{d,k}}{b_{d+p,k}}\,\frac{\E|P_\xi\El|^p}{\kappa_{k}^{p}}.
\end{equation}
\end{theorem}

Note that~\eqref{1124} is indeed a generalization of~\eqref{018} since $P_\xi\B^d=\B^k$ a.s. and $|\B^k|^p=\kappa_{k}^{p}$. For $k=1$ formula~\eqref{1124} was recently obtained in~\cite{lH14}.

For $p=1$, the right-hand side of~\eqref{1124} is proportional to the $k$-th intrinsic volume of $\El$ (see~\eqref{1309}), which  implies the following result (for details see Subsection~\ref{023}).

\begin{corollary}\label{024}
For  $X_0,\dots, X_{k}$ uniformly distributed in some non-degenerate ellipsoid $\El\subset\mathbb{R}^d$ we have
\[
\E|\conv(X_0,\dots,X_{k})|=\frac{1}{2^k}\,\frac{\left((d+1)!\right)^{k+1}}{\left((d+1)(k+1)\right)!}\,\left(\frac{\kappa_{d+1}^{k+1}}{\kappa_{(d+1)(k+1)}}\right)^2\,V_k(\El).
\]
\end{corollary}

Very recently, for  $X_0,\dots, X_{k}$ uniformly distributed in the unit ball $\B^d$, the formula for the distribution of $|\conv(X_0,\dots,X_k)|$ has been derived in~\cite{GKT17}. For a random variable $\eta$ and $\alpha_1,\alpha_2>0$ we write $\eta\sim B(\alpha_1,\alpha_2)$ to denote that $\eta$ has a Beta distribution with parameters $\alpha_1,\alpha_2$ and the density 
$$
\frac{\Gamma(\alpha_1+\alpha_2)}{\Gamma(\alpha_1)\,\Gamma(\alpha_2)}\,t^{\alpha_1-1}\,(1-t)^{\alpha_2-1},\quad t\in(0,1).
$$
It was shown  in~\cite{GKT17} that for  $X_0,\dots, X_{k}$ uniformly distributed in  $\B^d$,
\begin{equation}\label{005}
(k!)^2\,\eta(1-\eta)^k\,|\conv(X_0,\dots,X_k)|^2\eqdistr(1-\eta')^k\,\eta_1\cdots\eta_k,
\end{equation}
where $\eta,\eta',\eta_1,\dots,\eta_k$  are independent random variables independent of $X_0,\dots,X_k$ such that 
$$
\eta,\eta'\sim B\bigg(\frac d2+1,\frac{kd}{2}\bigg), \quad \eta_i\sim B\bigg(\frac{d-k+i}{2},\frac{k-i}{2}+1\bigg).
$$
Multiplying both sides of~\eqref{005} by $|P_\xi\El|^2/\kappa_k^2$ and applying Theorem~\ref{1305} and Corollary~\ref{009}  (for details see Subsection~\ref{020}) leads to the following generalization of~\eqref{005}.

\begin{theorem}\label{1716}
For $X_0,\dots, X_{k}$ uniformly distributed in some non-degenerate ellipsoid $\El\subset\mathbb{R}^d$ we have
\begin{align*}
(k!)^2\,\eta\,(1-\eta)^k\,|\conv(X_0,\dots,X_k)|^2&\eqdistr \kappa_k^{-2}\,(1-\eta')^k\,\eta_1\cdots\eta_k\,|P_\xi\El|^2\\
&\eqdistr(1-\eta')^k\,\eta_1\cdots\eta_k\,\left(\frac{\det\left(G_\lambda^\top G_\lambda\right)}{\det\left(G^\top G\right)}\right),
\end{align*}
where the matrices $G$ and $G_\lambda$ are defined in Corollary~\ref{009} and $\lambda_1,\dots,\lambda_d$ denote the length of semi-axes of $\El$.
\end{theorem}
Taking $k=1$ yields the distribution of the distance between two random points in $\El$.

\begin{corollary}
Under the assumptions of Theorem \ref{1716} we have
	\[
	\eta(1-\eta)\,\cdot|X_0-X_1|^2\eqdistr(1-\eta')\,\eta_1\,\left(\frac{\lambda_1^2N_1^2+\dots+\lambda_d^2N_d^2}{N_1^2+\dots+N_d^2}\right),
	\]
	where  $N_1,\dots,N_d$ are i.i.d. standard Gaussian variables.
\end{corollary}

\subsection{Integral geometry formulas}

Recall that $G_{d,k}$ and $A_{d,k}$ denote the linear and affine Grassmannians defined in Subsection~\ref{1733}.

For an arbitrary convex compact body $K$, $p>-d$, and $k=1$ it is possible to express~\eqref{1731} in terms of the lengths of the one-dimensional sections of $K$ (see~\cite{gC67} and~\cite{jK69}):
\[
\int\limits_{K^2}{|\bx_0-\bx_1|^p\,\dd\bx_0\dd\bx_1}=\frac{2d\kappa_d}{(d+p)\,(d+p+1)}\int\limits_{A_{d,1}}|K\cap E|^{p+d+1}\,\mu_{d,1}(dE).
\]  
This formula does not extend to $k>1$. The next theorem shows that for ellipsoids this is possible.
\begin{theorem}\label{2228}
For any  non-degenerate ellipsoid $\El\subset\mathbb{R}^d$, $k\in\{0,1,\dots,d\}$, and any real number $p>-d+k-1$ we have
\begin{multline}\label{1001}
\int\limits_{\El^{k+1}}{|\conv(\bx_0,\ldots,\bx_k)|^p\,\dd\bx_0\ldots \dd\bx_k}
\\=\frac{1}{(k!)^{p}}\frac{\kappa_{d+p}^{k+1}}{\kappa_k^{p+d+1}}\,\frac{\kappa_{k(d+p)+k}}{\kappa_{(k+1)(d+p)}}\,\frac{b_{d,k}}{b_{d+p,k}}\int\limits_{A_{d,k}}|\El\cap E|^{p+d+1}\,\mu_{d,k}(dE).
\end{multline}
\end{theorem}
The proof is given is Subsection~\ref{sec2228}.

Comparing this theorem with Theorem~\ref{2227} readily gives the following connection between the volumes of $k$-dimensional cross-sections and projections of ellipsoids.

\begin{theorem}\label{1115}
Under the  assumptions of Theorem~\ref{2228} we have
	\[
	\frac{\kappa_{d}^{k+1}}{\kappa_k^{d+1}}\,\frac{\kappa_{k(d+p)+k}}{\kappa_{k(d+p)+d}}\,\int\limits_{A_{d,k}}|\El\cap E|^{p+d+1}\,\mu_{d,k}(dE)=|\El|^{k+1}\,\int\limits_{G_{d,k}}|P_L\El|^p\,\nu_{d,k}(dL).
	\]
\end{theorem}

For $p=0$, we obtain the following integral formula.
\begin{corollary}
Under the  assumptions of Theorem~\ref{2228} we have
\begin{equation}\label{1113}
	\int\limits_{A_{d,k}}|\El\cap E|^{d+1}\mu_{d,k}(dE)=\frac{\kappa_k^{d+1}}{\kappa_{d}^{k+1}}\,\frac{\kappa_{d(k+1)}}{\kappa_{k(d+1)}}\,|\El|^{k+1}.
\end{equation}
\end{corollary}

This result may be regarded as an affine version of the following integral formula of Furstenberg and Tzkoni (see~\cite{FT71}):
$$
\int\limits_{G_{d,k}}|\El\cap L|^d\,\nu_{d,k}(dL)=\frac{\kappa_k^d}{\kappa_d^k}\,|\El|^k.
$$
Our next theorem generalizes this formula in the same way as~\eqref{1001} generalizes~\eqref{1113}. 

\begin{theorem}\label{1114}
	For any  non-degenerate ellipsoid $\El\subset\mathbb{R}^d$, $k\in\{0,1,\dots,d\}$, and any real number $p> -d+k$ we have
	\begin{equation}\label{2338}
	\int\limits_{\El^{k}}{|\conv(0,\bx_1,\ldots,\bx_k)|^p\,\dd\bx_1\ldots \dd\bx_k}
	=\frac{1}{(k!)^{p}}\frac{\kappa_{d+p}^{k}}{\kappa_k^{p+d}}\,\frac{b_{d,k}}{b_{d+p,k}}\,\int\limits_{G_{d,k}}|\El\cap L|^{p+d}\,\nu_{d,k}(dL).
	\end{equation}
\end{theorem}

In probabilistic language it may be formulated as
\[
\E|\conv(0,X_1,\dots,X_k)|^p=\frac{1}{(k!)^{p}}\frac{\kappa_{d+p}^{k}}{\kappa_k^{p+d}}\,\frac{b_{d,k}}{b_{d+p,k}}\,\E|\El\cap\xi|^{p+d},
\]
where $X_1,\dots,X_k$ are independent uniformly distributed random vectors in $\El$ and $\xi$ is a  uniformly chosen random $k$-dimensional linear subspace in $\R^d$.

\section{Proofs: Part~I}\label{2252}

\subsection{Proof of Theorem~\ref{1305} assuming Proposition~\ref{1843}}

Fist of all note that with probability one the equation 
	\[
	\frac{|P_\xi\El|}{\kappa_k}\cdot|\conv(X_0,\dots,X_{k})|=0
	\]
	holds if and only if
	\[
	|\conv(AX_0,\dots,AX_{k})|=0,
	\]
	which in turn is equivalent to
	\[
	\dim\conv(X_0,\dots, X_k)<k.
	\]
	Therefore to prove \eqref{1514} it is enough to show that the conditional distributions of $|\conv(AX_0,\dots,AX_{k})|$ and $\frac{|P_\xi\El|}{\kappa_k}\cdot|\conv(X_0,\dots,X_{k})|$ given $\dim\conv(X_0,\dots, X_k)=k$ are equal. Thus without loss of generality we can assume that the simplex  $\conv(X_0,\dots, X_k)$ is not degenerate with probability one:
\begin{equation}\label{2334}
	\dim\conv(X_0,\dots, X_k)=k\quad\text{a.s.}
\end{equation}
Our original proof was based on the Blaschke-Petkantschin formula and the characteristic function uniqueness theorem.\footnote{The original proof can be found in the first version of arXiv:1711.06578} Later, Youri Davydov has found a much simpler and nicer proof which also allowed to get rid of the assumption about the existence of the joint density of $X_0,\dots,X_{k}$. Let us present this proof.

Since the joint distribution of $X_0,\dots, X_{k}$ is spherically symmetric, we have for any orthogonal matrix $U$ 
\begin{align}\label{2043}
|\conv(AX_0,\dots,AX_{k})|&=|\conv(0, A(X_1-X_0),\dots,A(X_{k}-X_0)|\\
&\eqdistr|\conv(0, A(UX_1-UX_0),\dots,A(UX_{k}-UX_0)|.\notag
\end{align}
Now let $\Upsilon$ be a random orthogonal matrix chosen uniformly from $SO(n)$ with respect to the probabilistic Haar measure and independently of $X_0,\dots, X_{k}$. By~\eqref{2334}, with probability one the span of $X_1-X_0,\dots,X_{k}-X_0$ is a $k$-dimensional linear subspace of $\R^d$. Thus, the span
\[
\xi:=\SPAN(\Upsilon X_1-\Upsilon X_0,\dots,\Upsilon X_{k}-\Upsilon X_0)
\]
is a random uniformly chosen $k$-dimensional linear subspace in $\R^d$ independent of $X_0,\dots, X_{k}$. Applying Proposition~\ref{1843} to the vectors $\Upsilon X_1-\Upsilon X_0,\dots,\Upsilon X_{k}-\Upsilon X_0$ we obtain
\begin{align*}
|\conv(0, A(\Upsilon X_1&-\Upsilon X_0),\dots,A(\Upsilon X_{k}-\Upsilon X_0)|\\
&=\frac{|P_\xi\El|}{\kappa_k}\cdot|\conv(0,\Upsilon X_1-\Upsilon X_0,\dots,\Upsilon X_{k}-\Upsilon X_0)|\\
&=\frac{|P_\xi\El|}{\kappa_k}\cdot|\conv(\Upsilon X_0,\Upsilon X_1,\dots,\Upsilon X_{k})|\eqdistr\frac{|P_\xi\El|}{\kappa_k}\cdot|\conv(X_0, X_1,\dots,X_{k})|.
\end{align*}
Combining this with~\eqref{2043} for $U=\Upsilon$ finishes the proof.

\subsection{Proof of Proposition~\ref{1843}}
To avoid trivialities we assume that $\dim L=k$, i.e. $\bx_1,\ldots,\bx_k$ are in general position.
Let $\bp_1,\ldots,\bp_k\in\R^d$ be some orthonormal basis in $L$. Let $O_L$ and  $X$ denote $d\times k$ matrices  whose columns are $\bp_1,\ldots,\bp_k$ and $\bx_1,\dots,\bx_k$ respectively. It is easy to check that  $O_LO_L^{\top}$  is a $d\times d$ matrix corresponding to the orthogonal projection operator $P_L$. Thus,
\begin{equation}\label{1024}
O_LO_L^{\top}X=X.
\end{equation}

Recall that $\El$ is defined by~\eqref{2259}. It is known (see, e.g.,~\cite[Appendix~H]{Sch73}) that the orthogonal projection $ P_L\El$ is an ellipsoid in $L$ and
\begin{equation}\label{2304}
| P_L\El|=\kappa_k\Big[\det\left(O_L^{\top}HO_L\right)\Big]^{1/2},
\end{equation}
where 
\[
H:=A^\top A.
\]

A well-known formula for the volume of a $k$-dimensional arallelepiped implies that for any $\bx_1,\ldots,\bx_k\in\R^d$
\begin{equation}\label{1033}
|\conv(0,\bx_1,\ldots,\bx_k)|=\frac{1}{k!}\Big[\det\left(X^\top  X\right)\Big]^{1/2}.
\end{equation}
Therefore,  
\[
k!\,|\conv(0,A\bx_1,\ldots,A\bx_k)|=\Big[\det\left(\left(AX\right)^\top AX\right)\Big]^{1/2}=\Big[\det\left(X^\top HX\right)\Big]^{1/2}.
\]
Applying ~\eqref{1024} produces 
\begin{align*}
\det\left(X^\top HX\right)&=\det\left(X^\top O_LO_L^\top HO_LO_L^\top X\right)=\det\left(O_L^\top HO_L\right)\det\left(X^\top O_L\right)\det\left(O_L^\top X\right)\\
&=\det\left(O_L^\top HO_L\right)\det\left(X^\top O_LO_L^\top X\right)=\det\left(O_L^\top HO_L\right)\det\left(X^\top  X\right),
\end{align*}
which together with~\eqref{2304} and~\eqref{1033} finishes the proof.

\subsection{Proof of Corollary~\ref{009}}\label{016}

Denote by $G_1,\dots,G_k\in \R^d$ the columns of the matrix $G$. Hence, $AG_1,\dots, AG_k\in \R^d$ are the columns of the  matrix $AG$. Using Proposition \ref{1843} with $\bx_i=G_i$ and applying~\eqref{1033} to $G$ and $AG$ gives
\[
\Big[\det\big(G^\top A^\top A G\big)\Big]^{1/2}=\frac{|P_\eta\El|}{\kappa_k}\cdot\Big[\det\big(G^\top G\big)\Big]^{1/2},
\]
or
\[
\left({\frac{\det\big(G^\top A^\top AG\big)}{\det\big(G^\top G\big)}}\right)^{1/2}=\frac{|P_\eta\El|}{\kappa_k},
\]
where $\eta$ is the span of $G_1,\dots,G_k$. Since  $G_1,\dots,G_k$ are i.i.d. standard Gaussian vectors, $\eta$ is uniformly distributed in $G_{d,k}$ with respect to $\nu_{d,k}$ (given $\dim\eta=k$ which holds a.s.), therefore $\eta\eqdistr\xi$ and the corollary follows.

\subsection{Proofs of Theorem~\ref{2227} and Theorem~\ref{1716}}\label{020}
For any non-degenerate  ellipsoid $\El$ there exists a unique \emph{symmetric} positive-definite $d\times d$ matrix $A$ such that 
\[
\El=A\B^d=\left\{\bx\in\R^d:\|A^{-1}\bx\|\leq 1\right\}=\left\{\bx\in\R^d:{\bx^\top A^{-2}\bx}\leq 1\right\}.
\]
Since $X_0,\dots,X_{k}$ are i.i.d. random vectors uniformly distributed in $\El$, we have that $A^{-1}X_0,\dots,A^{-1}X_{k}$ are i.i.d. random vectors uniformly distributed in $\B^d$. It follows from Theorem~\ref{1305} that
\begin{multline}\label{2022}
|\conv(X_0,\dots,X_{k})|=|\conv(AA^{-1}X_0,\dots,AA^{-1}X_{k})|\\
\eqdistr|\conv(A^{-1}X_0,\dots,A^{-1}X_{k})|\,\frac{|P_\xi\El|}{\kappa_{k}}.
\end{multline}
Taking the $p$-th moment and applying~(\ref{018})  implies Theorem~\ref{2227}.

Now apply~\eqref{005} to $A^{-1}X_0,\dots,A^{-1}X_{k}$:
\[
(k!)^2\,\eta(1-\eta)^k\,|\conv(A^{-1}X_0,\dots,A^{-1}X_k)|^2\eqdistr(1-\eta')^k\,\eta_1\cdots\eta_k.
\]
Multiplying by $\frac{|P_\xi\El|}{\kappa_{k}^{p}}$ and applying~\eqref{2022} implies the first equation in Theorem~\ref{1716}. The second one follows from~\eqref{2306}.

\subsection{Proof of Corollary~\ref{024}}\label{023}

From Kubota's formula (see \eqref{1309}) and Theorem~\ref{2227} we have
\[
\E|\conv(X_0,\dots,X_{k})|=\alpha_{d,k}\,V_k(\El),
\]
where
\[
\alpha_{d,k}:=\frac{1}{k!}\,\frac{\kappa_{d+1}^{k+1}}{\kappa_d^{k+1}}\,\frac{\kappa_{k(d+1)+d}}{\kappa_{(k+1)(d+1)}}\,\frac{b_{d,k}}{b_{d+1,k}}\,\frac{\kappa_{d-k}}{{d \choose k}\,\kappa_{d}}.
\]
From the definition of $b_{d,k}$ (see (\ref{1846})) and $\kappa_p$ (see (\ref{1409})) we obtain
\begin{align*}
\alpha_{d,k} &= \frac{\kappa_{d+1}^{k+1}}{\kappa_d^{k+1}}\,\frac{\kappa_{k(d+1)+d}}{\kappa_{(k+1)(d+1)}}\,\frac{(d+1-k)!\,\kappa_{d-k+1}}{(d+1)!\,\kappa_{d+1}}\,\frac{\kappa_{d-k}}{\kappa_{d}}\\&= \frac{(d+1-k)!}{\pi^{k/2}(d+1)!}\,\left(\frac{\Gamma\left(\frac12d+1\right)}{\Gamma\left(\frac12(d+1)+1\right)}\right)^{k+1}
\\&\times\frac{\Gamma\left(\frac12(k+1)(d+1)+1\right)}{\Gamma\left(\frac12\left((k+1)d+k\right)+1\right)}\,\frac{\Gamma\left(\frac12(d+1)+1\right)}{\Gamma\left(\frac12(d-k+1)+1\right)}\,\frac{\Gamma\left(\frac12d+1\right)}{\Gamma\left(\frac12(d-k)+1\right)}.
\end{align*}
Using  Legendre's duplication formula for the Gamma function 
\[
\Gamma\left(z\right)\,\Gamma\left(z+\frac12\right)=2^{1-2z}\,\pi^{1/2}\,\Gamma\left(2z\right),
\]
the recursion $\Gamma\left(1+z\right)=z\,\Gamma\left(z\right)$, and the fact that $k,d\in\Z$ we obtain
\begin{align*}
\alpha_{d,k}&=\frac{(d-k)!}{\pi^{k/2}d!}\,\frac{\Gamma\left(\frac12(k+1)(d+1)+1\right)}{\Gamma\left(\frac12\left((k+1)d+k\right)+1\right)}\frac{\Gamma\left(\frac12d+\frac12\right)\,\Gamma\left(\frac12d+1\right)}{\Gamma\left(\frac12(d-k)+\frac12\right)\,\Gamma\left(\frac12(d-k)+1\right)}
\left(\frac{\Gamma\left(\frac12d+1\right)}{\Gamma\left(\frac12(d+1)+1\right)}\right)^{k+1}
\\&=\frac{1}{\left(2\sqrt{\pi}\right)^k}\,\frac{\Gamma\left(\frac12(k+1)(d+1)+1\right)}{\Gamma\left(\frac12\left((k+1)d+k\right)+1\right)}\,\left(\frac{\Gamma\left(\frac12d+1\right)}{\Gamma\left(\frac12(d+1)+1\right)}\right)^{k+1}\\
&=\frac{1}{\left(2\sqrt{\pi}\right)^k}\,\frac{\left(\Gamma\left(\frac12d+1\right)\,\Gamma\left(\frac12d+1+\frac12\right)\right)^{k+1}}{\Gamma\left(\frac12\left(kd+d+k\right)+1\right)\,\Gamma\left(\frac12(kd+k+d)+1+\frac12\right)}\,\left(\frac{\kappa_{d+1}^{k+1}}{\kappa_{(d+1)(k+1)}}\right)^2\\
&=\frac{1}{2^k}\,\frac{\left((d+1)!\right)^{k+1}}{\left((d+1)(k+1)\right)!}\,\left(\frac{\kappa_{d+1}^{k+1}}{\kappa_{(d+1)(k+1)}}\right)^2.
\end{align*}

\section{Proofs: Part~II}\label{037}

\subsection{Blaschke-Petkantschin formula}\label{2257}

In our further calculations we will need to integrate some non-negative measurable  function $h$ of $k$-tuples of points in $\R^d$. To this end, we integrate first over the $k$-tuples of points in a fixed $k$-dimensional linear subspace $L$ with respect to the product measure $\lambda_L^k$ and then integrate over $G_{d,k}$ with respect to  $\nu_{d,k}$. The corresponding transformation formula is known as the linear Blaschke-Petkantschin formula (see~ \cite[Theorem 7.2.1]{SW08}):
\begin{align}\label{eq3}
\int\limits_{(\R^d)^k}&{h(\bx_1,\ldots,\bx_k)\,\dd\bx_1\ldots \dd\bx_k}=\\&(k!)^{d-k}b_{d,k} \int\limits_{G_{d,k}}\int\limits_{L^k}h(\bx_1,\ldots,\bx_k)\,|\conv(0,\bx_1,\ldots,\bx_k)|^{d-k}\,\lambda_L(\dd\bx_1)\ldots\lambda_L(\dd\bx_k)\,\nu_{d,k}(\dd L),\notag
\end{align}
where $b_{d,k}$ is defined in~\eqref{1846}.

A similar affine version (see~\cite[Theorem 7.2.7]{SW08}) may be stated as follows:
\begin{align}\label{eq3_2}
\int\limits_{(\R^d)^{k+1}}&{h(\bx_0,\ldots,\bx_k)\,\dd\bx_0\ldots \dd\bx_k}=\\&(k!)^{d-k}b_{d,k} \int\limits_{A_{d,k}}\int\limits_{E^{k+1}}h(\bx_0,\ldots,\bx_k)\,|\conv(\bx_0,\ldots,\bx_k)|^{d-k}\,\lambda_E(\dd\bx_0)\ldots\lambda_E(\dd\bx_k)\,\mu_{d,k}(\dd E).\notag
\end{align}

\subsection{Proof of Theorem~\ref{2228}}\label{sec2228}
Let
\begin{multline*}
J:=\int\limits_{\El^{k+1}}{|\conv(\bx_0,\ldots,\bx_k)|^p\,\dd\bx_0\ldots \dd\bx_k}\\
=\int\limits_{(\R^d)^{k+1}}{|\conv(\bx_0,\ldots,\bx_k)|^{p}\,\prod\limits_{i=0}^k\ind_{\El}(\bx_i)\,\dd\bx_0\ldots \dd\bx_k}.
\end{multline*}

Using the affine Blaschke-Petkantschin formula (see~\eqref{eq3_2}) with
\[
h(\bx_0,\ldots,\bx_k):=|\conv(\bx_0,\ldots,\bx_k)|^{p}\,\prod\limits_{i=0}^k\ind_{\El}(\bx_i)
\]
yields
\begin{align*}
J&=(k!)^{d-k}b_{d,k}\,\int\limits_{A_{d,k}}\int\limits_{E^{k+1}}|\conv(\bx_0,\ldots,\bx_k)|^{p+d-k}\,\prod\limits_{i=0}^k\ind_{\El}(\bx_i)\,\lambda_E(\dd\bx_0)\ldots\lambda_E(\dd\bx_k)\,\mu_{d,k}(\dd E)\\
&=(k!)^{d-k}b_{d,k}\,\int\limits_{A_{d,k}}\int\limits_{(E\cap\El)^{k+1}}|\conv(\bx_0,\ldots,\bx_k)|^{p+d-k}\,\lambda_E(\dd\bx_0)\ldots\lambda_E(\dd\bx_k)\,\mu_{d,k}(\dd E).
\end{align*}
Now fix $E\in A_{d,k}$. Applying  Theorem~\ref{2227} to the ellipsoid $\El\cap E$ gives
\begin{multline*}
\frac{1}{|\El\cap E|^{k+1}}\,\int\limits_{(E\cap\El)^{k+1}}|\conv(\bx_0,\ldots,\bx_k)|^{p+d-k}\,\lambda_E(\dd\bx_0)\ldots\lambda_E(\dd\bx_k)\\
=\frac{1}{(k!)^{p+d-k}}\frac{\kappa_{d+p}^{k+1}}{\kappa_k^{k+1}}\,\frac{\kappa_{k(d+p)+k}}{\kappa_{(k+1)(d+p)}}\,\frac{b_{k,k}}{b_{d+p,k}}\,\frac{|\El\cap E|^{p+d-k}}{\kappa_k^{p+d-k}},
\end{multline*}
which leads to
\[
J=\frac{1}{(k!)^{p}}\frac{\kappa_{d+p}^{k+1}}{\kappa_k^{p+d+1}}\,\frac{\kappa_{k(d+p)+k}}{\kappa_{(k+1)(d+p)}}\,\frac{b_{d,k}}{b_{d+p,k}}\,\int\limits_{A_{d,k}}|\El\cap E|^{p+d+1}\,\mu_{d,k}(\dd E).
\]

\subsection{Proof of Theorem~\ref{1114}}\label{sec012}

The proof is  similar to the previous one. Let

\begin{multline*}
J:=\int\limits_{\El^{k}}{|\conv(0,\bx_1,\ldots,\bx_k)|^p\,\dd\bx_1\ldots \dd\bx_k}\\
=\int\limits_{(\R^d)^{k}}{|\conv(0,\bx_1,\ldots,\bx_k)|^{p}\,\prod\limits_{i=1}^k\ind_{\El}(\bx_i)\,\dd\bx_1\ldots \dd\bx_k}.
\end{multline*}

Using the linear Blaschke-Petkantschin formula (see~\eqref{eq3}) with
\[
h(\bx_1,\ldots,\bx_k):=|\conv(0,\bx_1,\ldots,\bx_k)|^{p}\,\prod\limits_{i=1}^k\ind_{\El}(\bx_i)
\]
gives
\begin{align}\label{91}
J&=(k!)^{d-k}b_{d,k}\,\int\limits_{G_{d,k}}\int\limits_{L^{k}}|\conv(0,\bx_1,\ldots,\bx_k)|^{p+d-k}\,\prod\limits_{i=1}^k\ind_{\El}(\bx_i)\,\lambda_L(\dd\bx_1)\ldots\lambda_L(\dd\bx_k)\,\nu_{d,k}(\dd L)\notag\\
&=(k!)^{d-k}b_{d,k}\,\int\limits_{G_{d,k}}\int\limits_{(L\cap\El)^{k}}|\conv(0,\bx_1,\ldots,\bx_k)|^{p+d-k}\,\lambda_L(\dd\bx_1)\ldots\lambda_L(\dd\bx_k)\,\nu_{d,k}(\dd L).
\end{align}
Fix $L\in G_{d,k}$. Since $\El\cap L$ is an ellipsoid,  there exists a linear transform $A_L:L\to\R^k$ such that $A_L(\El\cap L)=\B^k$. Applying the coordinate transform $\bx_i=A_L\by_i$, $i=1,2,\ldots,k$, we get
\begin{multline}\label{92}
\int\limits_{(L\cap\El)^{k}}|\conv(0,\bx_1,\ldots,\bx_k)|^{p+d-k}\,\lambda_L(\dd\bx_1)\ldots\lambda_L(\dd\bx_k)\\
=\frac{|\El\cap L|^{p+d}}{\kappa_k^{p+d}}\,\int\limits_{(\B^k)^{k}}|\conv(0,\by_1,\ldots,\by_k)|^{p+d-k}\,\dd\by_1\ldots\dd\by_k.
\end{multline}
It is known (see, e.g.,~\cite[Theorem~8.2.2]{SW08}) that
\begin{equation}\label{029}
\int\limits_{(\B^k)^{k}}|\conv(0,\by_1,\ldots,\by_k)|^{p+d-k}\,\dd\by_1\ldots\dd\by_k
=(k!)^{-p-d+k}\kappa_{d+p}^k\,\frac{b_{k,k}}{b_{d+p,k}}.
\end{equation}
Substituting~\eqref{029} and~\eqref{92} into~\eqref{91} finishes the proof.

\subsection{Acknowledgements}
The authors are grateful to Daniel Hug and G\"unter Last  for helpful discussions and suggestions which improved this paper. The authors also gratefully acknowledge Youri Davydov who considerably simplified the proof of Theorem~\ref{1305}.

The work was done with the financial support of the Bielefeld University (Germany). The work of F.G. and D.Z. is supported by the grant SFB 1283. The work of A.G. is supported by the grant IRTG 2235. The work of D.Z. is supported by the grant RFBR 16-01-00367 and by the Program of Fundamental Researches of Russian Academy of Sciences ``Modern Problems of Fundamental Mathematics''

\bibliographystyle{plainnat}

\end{document}